\documentclass[11pt,leqno]{article}
\usepackage{amsmath, amscd, amsthm, amssymb, graphics, xypic, mathrsfs,setspace,fancyhdr, times}
\usepackage[pagebackref=true]{hyperref}
\hypersetup{backref}

\newcommand{\tensor}{\otimes}
\newcommand{\colim}{\operatorname{colim}}
\newcommand{\Spec}{\operatorname{Spec}}
\newcommand{\isomto}{{\stackrel{\sim}{\;\longrightarrow\;}}}
\newcommand{\isomt}{{\stackrel{{\scriptscriptstyle{\sim}}}{\;\rightarrow\;}}}

\renewcommand{\O}{{\mathcal O}}
\renewcommand{\hom}{\operatorname{Hom}}

\newcommand{\cplx}{{\mathbb C}}
\newcommand{\Q}{{\mathbb Q}}
\newcommand{\Z}{{\mathbb Z}}
\newcommand{\aone}{{\mathbb A}^1}
\newcommand{\pone}{{\mathbb P}^1}

\newcommand{\gm}{{{\mathbb G}_{m}}}
\newcommand{\et}{\text{\'et}}
\newcommand{\ho}[1]{{\mathcal H}({#1})}

\newcommand{\dmeff}{{\mathbf{DM}}^{eff}_{k,-}}

\newcommand{\Shv}{{\mathcal Shv}}
\newcommand{\Sm}{{\mathcal Sm}}
\newcommand{\Cor}{{\mathcal Cor}}

\newcommand{\Mod}{{\mathcal Mod}}
\newcommand{\Ab}{{\mathcal Ab}}
\newcommand{\K}{{{\mathbf K}}}
\renewcommand{\H}{{{\mathbf H}}}

\newcommand{\hsnis}{{\mathcal H}_s((\Sm_k)_{Nis})}

\newcommand{\F}{{\mathcal F}}

\newcommand{\simpnis}{{\Delta}^{\circ}\Shv_{Nis}({\mathcal Sm}_k)}
\newcommand{\simpmod}{{\Delta}^{\circ}\Mod}

\newcounter{intro}
\theoremstyle{plain}
\newtheorem{thm}{Theorem}[section]

\newtheorem{lem}[thm]{Lemma}
\newtheorem{cor}[thm]{Corollary}
\newtheorem{prop}[thm]{Proposition}

\newtheorem*{thm*}{Theorem}
\newtheorem*{problem*}{Problem}

\newtheorem{thmintro}{Theorem}
\newtheorem{propintro}[thmintro]{Proposition}

\newtheorem{lemintro}[thmintro]{Lemma}

\theoremstyle{definition}
\newtheorem{defn}[thm]{Definition}

\newtheorem{notation}[thm]{Notation}

\theoremstyle{remark}
\newtheorem{rem}[thm]{Remark}
\newtheorem{remintro}[thmintro]{Remark}

\newtheorem{ex}[thm]{Example}

\numberwithin{equation}{section}

\begin{document}
\pagestyle{fancy}
\renewcommand{\sectionmark}[1]{\markright{\thesection\ #1}}
\fancyhead{}
\fancyhead[LO,R]{\bfseries\footnotesize\thepage}
\fancyhead[LE]{\bfseries\footnotesize\rightmark}
\fancyhead[RO]{\bfseries\footnotesize\rightmark}
\chead[]{}
\cfoot[]{}
\setlength{\headheight}{1cm}

\author{\begin{small}Aravind Asok\thanks{Aravind Asok was partially supported by National Science Foundation Awards DMS-0900813 and DMS-0966589.}\end{small} \\ \begin{footnotesize}Department of Mathematics\end{footnotesize} \\ \begin{footnotesize}University of Southern California\end{footnotesize} \\ \begin{footnotesize}Los Angeles, CA 90089-2532 \end{footnotesize} \\ \begin{footnotesize}\url{asok@usc.edu}\end{footnotesize}             }

\title{{\bf Birational invariants and $\aone$-connectedness}}
\date{}
\maketitle

\begin{abstract}
We study some aspects of the relationship between ${\mathbb A}^1$-homotopy theory and birational geometry.  We study the so-called ${\mathbb A}^1$-singular chain complex and zeroth ${\mathbb A}^1$-homology {\em sheaf} of smooth algebraic varieties over a field $k$.  We exhibit some ways in which these objects are similar to their counterparts in classical topology and similar to their motivic counterparts (the (Voevodsky) motive and zeroth Suslin homology sheaf).  We show that if $k$ is infinite the zeroth ${\mathbb A}^1$-homology sheaf is a birational invariant of smooth proper varieties, and we explain how these sheaves control various cohomological invariants, e.g., unramified \'etale cohomology.  In particular, we deduce a number of vanishing results for cohomology of ${\mathbb A}^1$-connected varieties.  Finally, we give a partial converse to these vanishing statements by giving a characterization of ${\mathbb A}^1$-connectedness by means of vanishing of unramified invariants.
\end{abstract}

\begin{footnotesize}
\tableofcontents
\end{footnotesize}

\section{Introduction}
In this paper, we continue to investigate the relationship between birational geometry and connectedness in the sense of $\aone$-homotopy theory that was initiated in \cite{AM}.  Developing some ideas of \cite[\S 4]{AM}, we study cohomological consequences of homotopical connectivity hypotheses and, more specifically, vanishing results for various types of cohomological invariants such as unramified \'etale cohomology.  For the most part, however, the results of this paper are logically independent of \cite{AM}.

In the context of the Morel-Voevodsky $\aone$-homotopy theory of smooth schemes over a field $k$ \cite{MV}, one may associate with any smooth scheme $X$ a Nisnevich sheaf of sets, denoted $\pi_0^{\aone}(X)$, called the sheaf of $\aone$-connected components of $X$.  A smooth scheme $X$ is called {\em $\aone$-connected} if $\pi_0^{\aone}(X)$ is isomorphic to the constant $1$-point sheaf (see Definition \ref{defn:aoneconnected} for more details).  In \cite[Definition 2.2.2 and Corollary 2.4.4]{AM} it was shown that varieties that are $\aone$-connected are ``nearly rational in a strong sense."  For example, if $k$ has characteristic $0$, then stably $k$-rational smooth proper varieties or, more generally, Saltman's retract $k$-rational smooth proper varieties are $\aone$-connected.

Many schemes of interest are not $\aone$-connected and providing an explicit description of $\pi_0^{\aone}(X)$ for any arbitrary smooth variety seems difficult at the moment.  For that reason, we seek to understand auxiliary invariants that are controlled by the sheaf of $\aone$-connected components.  Here, we study the homological counterpart of $\pi_0^{\aone}(X)$: the zeroth $\aone$-homology sheaf, denoted $\H_0^{\aone}(X)$ (see Definition \ref{defn:aonehomology}).  While in some ways $\aone$-homology is similar to the more familiar Suslin homology (see, e.g., \cite{SuslinVoevodsky}), the two theories are different.  While the zeroth Suslin homology sheaf (see Definition \ref{defn:suslinhomology}) of a smooth proper $k$-scheme $X$ can be described using the Chow group of $0$-cycles on $X$, as we explain, the zeroth $\aone$-homology is more closely related to $R$-equivalence classes in $X$.

Intuitively speaking, the sheaf $\pi_0^{\aone}(X)$ formalizes the idea of algebraic path components of $X$, where algebraic paths are interpreted as chains of affine lines.  Indeed \cite[Theorem 6.2.1]{AM} proves that if $X$ is smooth and proper over $k$, then for any finitely generated separable extension $L/k$, the set of sections $\pi_0^{\aone}(X)(L)$ can be identified with the set of $R$-equivalence classes in $X(L)$ in the sense of Manin.  While this result does not identify the whole sheaf $\pi_0^{\aone}(X)$, it suggests that $\pi_0^{\aone}(X)$ is a birational invariant of smooth proper $k$-varieties.

By analogy with topology, one might expect that the sheaf $\H_0^{\aone}(X)$ should be the free abelian group on the $\aone$-connected components of $X$.  Mirroring the expected behavior of $\pi_0^{\aone}(X)$ above, one might also expect that $\H_0^{\aone}(X)$ is a birational invariant for smooth proper $k$-varieties.  Both of these statements are true provided that in the first statement one interprets the expression ``free abelian group" correctly, and in the second statement one restricts $k$ appropriately (the corresponding statement for Suslin homology is well known).  Precisely, we prove the following results.

\begin{propintro}[See Proposition \ref{prop:dependsonpi0}]
\label{propintro:universalproperty}
Suppose $k$ is a field and $X$ is a smooth $k$-scheme.  The morphism $\H_0^{\aone}(X) \to \H_0^{\aone}(\pi_0^{\aone}(X))$, induced by the canonical morphism $X \to \pi_0^{\aone}(X)$, is an isomorphism of Nisnevich sheaves of abelian groups.
\end{propintro}

\begin{remintro}
By Lemma \ref{lem:h0universal}, $\H_0^{\aone}(X)$ is the free strictly $\aone$-invariant sheaf (see Definition \ref{defn:strictaoneinvariance}) of abelian groups generated by $X$.  Thus, this theorem says that $\H_0^{\aone}(X)$ is the free strictly $\aone$-invariant sheaf of groups on $\pi_0^{\aone}(X)$.  One particularly useful consequence of this result is that if the morphism $\H_0^{\aone}(X) \to \H_0^{\aone}(\Spec k)$ is not an isomorphism, then $X$ is $\aone$-disconnected.
\end{remintro}

\begin{thmintro}[See Theorem \ref{thm:birationalclass}]
\label{thmintro:birationalinvariance}
Suppose $k$ is an infinite field.  If $X$ and $X'$ are stably $k$-birationally equivalent smooth proper schemes, then $\H_0^{\aone}(X) \cong \H_0^{\aone}(X')$.
\end{thmintro}

For any field $k$, the sheaf $\H_0^{\aone}(\Spec k)$ is isomorphic to $\Z$ (see Example \ref{ex:aonehomologyofapoint}) and thus coincides with the zeroth Suslin homology of a point.  For a general smooth $k$-scheme $X$, various classical stable birational invariants are related to $\H_0^{\aone}(X)$.  Suppose $n$ is an integer coprime to the characteristic of $k$, and $L/k$ is a finitely generated separable field extension.  Consider the functor on $k$-algebras defined by $A \mapsto H^i_{\et}(A,\mu_n^{\tensor j})$ (we abuse terminology and write $A$ instead of $\Spec A$ for notational convenience).  Given a discrete valuation $\nu$ of $L/k$ with associated valuation ring $A$, one says that a class $\alpha \in H^i_{\et}(L,\mu_n^{\tensor j})$ is unramified at $\nu$ if $\alpha$ lies in the image of the restriction map $H^i_{\et}(A,\mu_n^{\tensor j}) \to H^i_{\et}(L,\mu_n^{\tensor j})$.  For any integers $i,j$, Colliot-Th\'el\`ene and Ojanguren \cite{CTO} define the unramified cohomology group $H^i_{ur}(L/k,\mu_n^{\tensor j})$ as the subgroup of $H^i_{\et}(L,\mu_n^{\tensor j})$ consisting of those classes $\alpha$ that are unramified at every discrete valuation of $L$ trivial on $k$.  Colliot-Th\'el\`ene and Ojanguren also proved that the groups $H^i_{ur}(L/k,\mu_n^{\tensor j})$ are stable birational invariants of smooth proper varieties (see \cite[Proposition 1.2]{CTO}).  For another point of view on these statements see \cite[Theorem 4.1.1]{CTPurity}.  The next result demonstrates that the groups $H^i_{ur}(L/k,\mu_n^{\tensor j})$ are controlled by $\H_0^{\aone}(X)$; in essence this was observed by Gabber ({\em cf.} \cite[Remark 1.1.3]{CTO}).

\begin{lemintro}[See Lemma \ref{lem:unramifiedrelation}]
\label{lemintro:unramifiedetalecohomology}
Suppose $k$ is a field, and $n$ is an integer coprime to the characteristic of $k$.  If $X$ is a smooth proper $k$-scheme, then there is a canonical bijection
\[
H^i_{ur}(X,\mu_n^{\tensor j}) \isomto \hom(\H_0^{\aone}(X),{\mathbf H}^i_{\et}(\mu_n^{\tensor j})),
\]
where the group on the right hand side is computed in the category of Nisnevich sheaves of abelian groups, and the sheaf ${\mathbf H}^i_{\et}(\mu_n^{\tensor j})$ is defined in \textup{Example \ref{ex:unramifiedetalecohomology}}.
\end{lemintro}

Proposition \ref{propintro:universalproperty} shows that if $X$ is $\aone$-connected, then $\H_0^{\aone}(X)$ is trivial.  Thus, for example, non-triviality of unramified \'etale cohomology can be used to detect $\aone$-disconnectedness.  More generally, via Lemma \ref{lem:h0universal}, we will see that $\H_0^{\aone}(X)$ is a ``universal unramified invariant" for smooth proper schemes, in an appropriate sense; see Lemma \ref{lem:unramifiedinvariants} for a precise statement.  For example, $\H_0^{\aone}(X)$ controls unramified Milnor K-theory (see Example \ref{ex:unramifiedMilnorktheory}) and the unramified Witt sheaf (see Example \ref{ex:unramifiedwittgroup}); this point of view is developed in \S \ref{s:unramifiedelements}.

All of the theories described in the previous paragraph have transfers of an appropriate kind, and the first two are even controlled by Suslin homology. However, the zeroth $\aone$-homology sheaf controls unramified invariants that do not possess transfers.  Using this additional information, if $X$ is smooth and proper we can show that $\aone$-connectedness is characterized by vanishing of unramified invariants.  More precisely, we prove the following result.

\begin{thmintro}[See Theorem \ref{thm:characterization}]
\label{thmintro:characterization}
If $k$ is a field, and $X$ is a smooth proper $k$-scheme, then $X$ is $\aone$-connected if and only if the canonical morphism $\H_0^{\aone}(X) \to \Z$ is an isomorphism; the same statement holds with rational coefficients.
\end{thmintro}

The use of strictly $\aone$-invariant sheaves that do not possess transfers seems essential.  Indeed, Suslin homology (even with integral coefficients) cannot detect $\aone$-connectedness, in large part because of its inability to see the difference between rational points and $0$-cycles of degree $1$.  We recall an example of Parimala (see Example \ref{ex:parimala}), pointed out to us by Sasha Merkurjev, showing that even if $X$ is a smooth projective variety such that the degree morphism $\H_0^S(X) \to \Z$ is an isomorphism, $X$ need not be $\aone$-connected.

Section \ref{s:aonehomology} is devoted to briefly reviewing aspects of $\aone$-homotopy theory, Voevodsky's theory of motives, and $\aone$-homology theory; in particular, we fix our notation for the rest of the paper.  Section \ref{s:birationalproperties} studies the birational properties of the zeroth $\aone$-homology and Suslin homology sheaves and relates these two objects to the zeroth $\aone$-homotopy sheaf.  Finally, Section \ref{s:unramifiedelements} provides a field theoretic point of view useful for studying $\aone$-homology and Suslin homology sheaves together with the universality statement alluded to above.

\subsubsection*{Relationship with other work}
This work is part of a sequence of papers including \cite{AH,AH2,ABK} studying the $\aone$-homology sheaf, its relationship with rational points, and rationality questions.  In \cite{AH}, we prove that if $X$ is smooth and proper over a field, then $\H_0^{\aone}(X)$ detects rational points.  More precisely \cite[Corollary 2.9]{AH} states that a smooth proper variety $X$ has a $k$-rational point if and only if the canonical map $\H_0^{\aone}(X) \to \Z$ is an epimorphism.  On the other hand $\pi_0^{\aone}(X)$ controls the unramified cohomology of smooth varieties $X$.  In \cite{ABK}, we produce rationally connected, non-rational smooth proper varieties $X$ where non-rationality is detected by a degree $n$ unramified cohomology class, but cannot be detected by lower degree invariants; the connection with this work is mentioned at the end of Section \ref{s:unramifiedelements}.

\subsubsection*{General conventions}
Throughout the paper $k$ will be a fixed base field.  Write $\Sm_k$ for the category of schemes that are smooth, separated and have finite type over $k$.  When we use the word {\em sheaf} without modification, we always mean Nisnevich sheaf of sets on $\Sm_k$.

\subsubsection*{Acknowledgements}
We were led to the study of cohomological vanishing properties of $\aone$-connected schemes by B. Bhatt's observation that the Brauer group of a $\aone$-connected smooth scheme is trivial (this is unpublished, but see \cite[Theorem 4.3]{Gille}), which suggested that higher cohomological invariants could obstruct $\aone$-connectedness as well.  This work also owes an intellectual debt to Fabien Morel, who has stressed the importance of $\aone$-homological invariants; we thank him for much encouragement and many discussions through the course of this work.  Finally, we thank Christian H\"asemeyer for many discussions around this topic, Sasha Merkurjev for pointing out Example \ref{ex:parimala}, and Jean-Louis Colliot-Th\'el\`ene for helpful comments and correspondence.

%\section{Birational invariants and $\aone$-homological algebra}
\section{$\aone$-homotopy, $\aone$-homology and Suslin homology (sheaves)}
\label{s:aonehomology}
We review the construction and basic properties of $\aone$-derived categories and $\aone$-homology as sketched in \cite{MStable} and developed in \cite{MField}.  For completeness, we also give a detailed comparison between $\aone$-homology and Suslin homology sheaves, which was alluded to in \cite{MICM} but has not been developed in the literature (in detail).  For more discussion of the homological algebra underlying the $\aone$-derived category, we refer the reader to \cite[\S 4]{CisinskiDeglise1}.  The results stated in this section are essential to the formulation and proofs of results in subsequent sections.

\subsubsection*{Simplicial homotopy categories}
Let $\simpnis$ denote the category of simplicial Nisnevich sheaves on $\Sm_k$; we will refer to objects in this category as $k$-spaces, or simply as spaces if $k$ is clear from context.  The Yoneda embedding provides a fully-faithful functor $\Sm_k \to \simpnis$.  We use this to identify $\Sm_k$ with a full subcategory of $\simpnis$, and systematically abuse notation by denoting a smooth scheme and the corresponding simplicial sheaf (the sheaf of $n$-simplices is the Nisnevich sheaf represented by the scheme, and all face and degeneracy morphisms are the identity morphism) by the same roman letter.  Generally, we use calligraphic letters (e.g., ${\mathcal X},{\mathcal Y}$) for objects of $\simpnis$.

The category $\simpnis$ admits a proper closed model structure where the cofibrations are monomorphisms, the weak equivalences are those morphisms of simplicial sheaves that stalkwise induce weak equivalences of the corresponding simplicial sets, and the fibrations are those morphisms having the right lifting property with respect to morphisms that are simultaneously cofibrations and weak equivalences (see, e.g., \cite[\S2 Theorem 1.4]{MV}).  The resulting model structure is called the injective local model structure, or the Joyal-Jardine model structure.  The {\em simplicial homotopy category}, denoted $\hsnis$, is the homotopy category of this model structure.  Throughout, we write $[{\mathcal X},{\mathcal Y}]_s$ for $\hom_{\hsnis}({\mathcal X},{\mathcal Y})$.

\subsubsection*{Derived categories of sheaves of $R$-modules}
For a commutative unital ring $R$, we let $\Mod_k(R)$ denote the category of Nisnevich sheaves of $R$-modules.  Similarly, we let $\simpmod_k(R)$ denote the category of simplicial Nisnevich sheaves of $R$-modules.  Given any object ${\mathcal X} \in \simpnis$, write $R({\mathcal X})$ for the Nisnevich sheaf of $R$-modules freely generated by the simplices of ${\mathcal X}$; $R({\mathcal X})$ an object of $\simpmod_k(R)$.  This construction defines a functor $R(\cdot): \simpnis \to \simpmod_k(R)$ that is left adjoint to the forgetful functor $\simpmod_k(R) \to \simpnis$.

Let $Ch_{\geq 0}(\Mod_k(R))$ denote the category of chain complexes (differential of degree $-1$) of Nisnevich sheaves of $R$-modules situated in degrees $\geq 0$.  There is a functor of normalized chain complex $N(\cdot): \simpmod_k(R) \to Ch_{\geq 0}(\Mod_k(R))$. The sheaf theoretic Dold-Kan correspondence produces an adjoint equivalence $K(\cdot): Ch_{\geq 0}(\Mod_k(R)) \to \simpmod_k(R)$.

Let $Ch_{-}(\Mod_k(R))$ denote the category of bounded below chain complexes of Nisnevich sheaves of $R$-modules; objects in this category will be referred to simply as {\em complexes}.  The category $\Sm_k$ is essentially small, so the category $Ch_{-}(\Mod_k(R))$ is the category of bounded below complexes in a Grothendieck abelian category.  Therefore, results of Beke imply that $Ch_{-}(\Mod_k(R))$ can be equipped with a model category structure where cofibrations are monomorphisms, weak equivalences are quasi-isomorphisms, and fibrations are those morphisms having the right lifting property with respect to morphisms that are simultaneously cofibrations and weak equivalences (see \cite[Proposition 3.13]{Beke}).  This model structure---the injective local model structure---has homotopy category the bounded below derived category $D_-(\Mod_k(R))$.  We denote by $((-)^f,\theta)$ a fixed fibrant resolution functor, i.e., $(-)^f$ is an endofunctor of $Ch_{-}(\Mod_k(R))$ and $\theta: Id \to (-)^f$ is a natural transformation such that if $A$ is complex, the induced map $A \to A^f$ is a quasi-isomorphism and monomorphism and $A^f$ is a fibrant complex.  The homotopy category $D_-(\Mod_k(R))$ is a triangulated category with the usual shift functor.

\begin{notation}
We use {\em homological conventions} for complexes.  More precisely, if $C_*$ is a complex, then the shift functor satisfies $C_*[1] = C_{*+1}$ so that $H_i(C_*[1]) = H_{i-1}(C_*)$; this convention will be justified in the next subsection.  Any complex $C_*$ can be considered as a cohomological complex $C^*$ with $C^i = C_{-i}$; we use this convention when computing hypercohomology.
\end{notation}

\subsubsection*{Hurewicz theory}
If ${\mathcal X}$ is a space, set $C_*({\mathcal X},R) := N((R({\mathcal X})))$.  The assignment ${\mathcal X} \to C_*({\mathcal X},R)$ provides a functor
\[
\simpnis \longrightarrow Ch_{\geq 0}(\Mod_k(R)).
\]
This functor sends monomorphisms to monomorphisms, and sends weak equivalences to quasi-isomorphisms.  Thus, it descends to a functor
\[
\hsnis \longrightarrow D_-(\Mod_k(R)).
\]
There is a corresponding version of this functor in the setting of pointed spaces as well.  If ${\mathcal X}$ is a space, the structure morphism ${\mathcal X} \to \Spec k$ induces a morphism of complexes $C_*({\mathcal X},R) \to C_*(\Spec k,R)$; we let $\tilde{C}_*({\mathcal X},R)$ denote the kernel of this morphism.  If $\mathcal{X}$ is pointed, then $C_*({\mathcal X},R) \to C_*(\Spec k,R)$ is split, and $\tilde{C}({\mathcal X},R)$ is a summand of $C_*({\mathcal X},R)$.

In the other direction, the adjoint $K(\cdot)$ (coming from the Dold-Kan correspondence) composed with the inclusion $\simpmod_k(R) \to \simpnis$ produces a functor
\[
Ch_{\geq 0}(\Mod_k(R)) \longrightarrow \simpnis.
\]
This composite functor sends quasi-isomorphisms to weak equivalences, and using properties of adjunctions can be shown to preserve fibrations as well.  In fact, there is an adjunction
\begin{equation}
\label{eqn:simplicialderivedadjunction}
[{\mathcal X},K(A,n)]_s \isomto \hom_{D_-(\Mod_k(R))}(R({\mathcal X}),A[n]),
\end{equation}
which we use freely in the sequel.

We set $H_i({\mathcal X},R) := H_i(C_*({\mathcal X},R))$ and if ${\mathcal X}$ is pointed $\tilde{H}_i({\mathcal X},R) := H_i(\tilde{C}_*({\mathcal X},R))$.  If $S^1_s$ denotes the constant sheaf defined by the simplicial circle, we let $\Sigma^1_s {\mathcal X} = \Sigma^1_s \wedge {\mathcal X}$.  It is not hard to check that $\tilde{H}_i(\Sigma^1_s {\mathcal X},R) = \tilde{H}_{i-1}({\mathcal X},R)$.

\subsubsection*{$\aone$-homotopy categories}
The $\aone$-homotopy category, constructed in \cite[\S 2 Theorem 3.2]{MV}, is obtained as a categorical localization of $\simpnis$.  Recall that a space ${\mathcal X}$ is called {\em $\aone$-local} if for any space ${\mathcal Y}$ the canonical map
\[
[{\mathcal Y},{\mathcal X}]_{s} \longrightarrow [{\mathcal Y} \times \aone,{\mathcal X}]_s
\]
is a bijection.  A morphism $f: {\mathcal X} \to {\mathcal Y}$ is an {\em $\aone$-weak equivalence} if the induced map $[{\mathcal Y},{\mathcal Z}]_{s} \to [{\mathcal X},{\mathcal Z}]_s$ is a bijection for all $\aone$-local spaces $\mathcal{Z}$.  The category $\simpnis$ can be equipped with a model structure where weak equivalences are $\aone$-weak equivalences, cofibrations are monomorphisms and fibrations are those morphisms having the right lifting property with respect to morphisms that are simultaneously cofibrations and $\aone$-weak equivalences.  We write $\ho{k}$ for the resulting homotopy category, and $[{\mathcal X},{\mathcal Y}]_{\aone}$ for $\hom_{\ho{k}}({\mathcal X},{\mathcal Y})$.  The full subcategory of $\hsnis$ spanned by $\aone$-local objects can be taken as a model for the $\aone$-homotopy category.  Thus, if ${\mathcal X}$ is $\aone$-local, we have $[{\mathcal Y},{\mathcal X}]_s = [{\mathcal Y},{\mathcal X}]_{\aone}$; we use this freely in the sequel.

\begin{defn}
\label{defn:aoneconnected}
Suppose ${\mathcal X}$ is a space. The {\em sheaf of $\aone$-connected components} of $\mathcal{X}$, denoted $\pi_0^{\aone}({\mathcal X})$ is the Nisnevich sheaf associated with the presheaf $U \mapsto [U,{\mathcal X}]_{\aone}$.  A space ${\mathcal X}$ is {\em $\aone$-connected} if the canonical morphism $\pi_0^{\aone}({\mathcal X}) \to \ast$ (where $\ast$ is the constant $1$ point sheaf) is an isomorphism.
\end{defn}

\subsubsection*{$\aone$-derived categories}
There is an analogous abelianized version of the $\aone$-homotopy category; we recall the basic definitions, which were originally introduced by Morel.

\begin{defn}
\label{defn:strictaoneinvariance}
A complex $A$ of Nisnevich sheaves of abelian groups is {\em $\aone$-local} if for any smooth $k$-scheme $U$, and every integer $n$
\[
{\mathbb H}^n_{Nis}(U,A) \to {\mathbb H}^n_{Nis}(U \times \aone,A)
\]
is an isomorphism.   If $A$ is simply a sheaf of abelian groups, say $A$ is {\em strictly $\aone$-invariant} if it is $\aone$-local viewed as a complex of sheaves.
\end{defn}

\begin{rem}
\label{rem:stronginvariance}
A sheaf of sets ${\mathcal S}$ is {\em $\aone$-invariant} if for any smooth scheme $U$, the map ${\mathcal S}(U) \to {\mathcal S}(U \times \aone)$ is a bijection.  A sheaf of groups ${\mathcal G}$ is {\em strongly $\aone$-invariant} if for any smooth scheme $U$, and any integer $i \in \{0,1\}$, the maps $H^i(U,{\mathcal G}) \to H^i(U \times \aone,{\mathcal G})$ are bijections.
\end{rem}

Denote the localizing subcategory of $D_-(\Mod_k(R))$ generated by complexes of the form $R(X \times \aone) \to R(X)$ for smooth schemes $X$ by $T(\aone,R)$.  By the theory of localizing categories, a quotient category $D_-(\Mod_k(R))/T(\aone,R)$ exists; this category, denoted $D_{\aone}(k,R)$, is called the $\aone$-derived category.  Let $D_-(\Mod_k(R))^{\aone-loc} \subset D_-(\Mod_k(R))$ denote the full-subcategory consisting of $\aone$-local complexes.  This inclusion admits a left adjoint
\[
L_{\aone}: D_-(\Mod_k(R)) \longrightarrow D_-(\Mod_k(R))^{\aone-loc}
\]
that can be used to identify, up to equivalence, $D_-(\Mod_k(R))/T(\aone,R)$ with $D_-(\Mod_k(R))^{\aone-loc}$.  The functor $L_{\aone}$ is called the $\aone$-localization functor; for more details see, e.g., \cite[Proposition 4.3]{CisinskiDeglise1}.  For simplicity, we will write $D_{\aone}(k)$ for $D_{\aone}(k,\Z)$.  Observe that if $A$ is an $\aone$-local complex, then $K(A,n)$ is an $\aone$-local space by the adjunction of \ref{eqn:simplicialderivedadjunction}.

\begin{defn}
\label{defn:aonehomology}
The {\em $\aone$-singular chain complex of ${\mathcal X}$ with $R$-coefficients}, denoted $C_*^{\aone}({\mathcal X},R)$, is defined to be the $\aone$-localization $L_{\aone}(C_*({\mathcal X},R))$.  We write $C_*^{\aone}({\mathcal X})$ for $C_*^{\aone}({\mathcal X},\Z)$.  The {\em $\aone$-homology sheaves of ${\mathcal X}$ with $R$-coefficients} are defined by $\H_i^{\aone}({\mathcal X},R) := H_i(C_*^{\aone}({\mathcal X},R))$.  If ${\mathcal X}$ is pointed, the {\em reduced $\aone$-homology sheaves of ${\mathcal X}$ with $R$-coefficients} are defined by $\tilde{\H}_i^{\aone}({\mathcal X},R) := H_i(L_{\aone}(\tilde{C}_*({\mathcal X},R)))$.
\end{defn}

\begin{ex}
\label{ex:aonehomologyofapoint}
Suppose $k$ is a field.  The complex $C_*(\Spec k)$ is just the Nisnevich sheaf $\Z$ placed in degree $0$.  The Nisnevich sheaf $\Z$ is even Nisnevich flasque so all higher Nisnevich cohomology of a smooth $k$-scheme with coefficients in $\Z$ vanishes.  Since the functor $\Z(\cdot)$ is clearly $\aone$-invariant, it follows immediately that $\Z$ is strictly $\aone$-invariant.  As a consequence $\H_0^{\aone}(\Spec k) = \Z$ and all the higher $\aone$-homology sheaves of $\Spec k$ vanish, as topological intuition suggests.
\end{ex}

\begin{rem}[Mayer-Vietoris]
If $X$ is a smooth scheme, there is a Nisnevich Mayer-Vietoris sequence allowing computationg of $\H_i^{\aone}(X)$ from pieces of an open cover. This follows immediately from \cite[\S Remark 1.7]{MV} together with the fact that the $\aone$-localization functor is exact.
\end{rem}

\subsubsection*{Sheaves with transfers and Suslin homology of motives}
Again, let $R$ be a commutative unital ring.  Write $Cor_k(X,Y)$ for the free $R$-module generated by integral closed subschemes of $X \times Y$ that are finite and surjective over a component of $X$ (an element of this group is referred to as a {\em finite correspondence}).  We let $R_{tr}(X)$ denote the presheaf on $\Sm_k$ defined by $U \mapsto Cor_k(U,X)$.  If $X$ is a smooth scheme, the functor $U \mapsto R_{tr}(X)(U \times \Delta^{\bullet})$ defines a simplicial abelian group for which we write $C_*^S(R_{tr}(X))$.  By definition $C_*^S(R_{tr}(X))$ is situated in positive (homological) degrees.

Write $\Cor_k$ for the category whose objects are smooth schemes, and morphisms are finite correspondences from between smooth schemes (this category is described in detail in \cite[Chapter 1]{MVW}).  There is a functor $\Sm_k \to \Cor_k$ that sends an element of $\hom(X,Y)$ to the finite correspondence defined by the graph.  An additive contravariant functor from $\Cor_k$ to abelian groups is called a presheaf with transfers.  Any presheaf with transfers can be viewed as a presheaf of abelian groups on $\Sm_k$ by restriction to $\Sm_k$.  A presheaf with transfers whose restriction to $\Sm_k$ is a Nisnevich sheaf is called a Nisnevich sheaf with transfers.  The presheaves $R_{tr}(X)$ are all Nisnevich sheaves with transfers (\cite[Lemma 6.2]{MVW}).

One can define a derived category of Nisnevich sheaves with transfers: take the homotopy category of complexes of Nisnevich sheaves with transfers and localize at the quasi-isomorphisms.  A complex of Nisnevich sheaves with transfers is $\aone$-local if it is $\aone$-local as a complex of Nisnevich sheaves after forgetting the transfers.  A {\em strictly $\aone$-invariant sheaf with transfers} is a sheaf with transfers that is $\aone$-local when viewed as a complex of Nisnevich sheaves with transfers.  Taking the quotient of the derived category of Nisnevich sheaves with transfers by the localizing subcategory generated by $\aone$-local complexes of Nisnevich sheaves with transfers one obtains Voevodsky's (``big") derived category of motives $\dmeff$.  We refer the reader to \cite[\S 13]{MVW} for a much more detailed discussion of this construction.  The techniques given so far are sufficient to define Suslin homology for a smooth scheme, but we will need to define Suslin homology of an arbitrary space for some later statements.  To do this, recall the following result.

\begin{lem}[{\cite[\S 2 Lemma 1.16]{MV}}]
There is a pair $(\Phi_{rep},\theta)$ consisting of an endofunctor $\Phi_{rep}: \simpnis \to \simpnis$ and a natural transformation $\theta: \Phi_{rep} \to Id$ such that for any simplicial sheaf ${\mathcal X}$, $\Phi_{rep}({\mathcal X})_n$ is a coproduct of representable sheaves, and $\Phi_{rep}({\mathcal X}) \to {\mathcal X}$ is a simplicial weak equivalence and stalkwise a fibration of simplicial sets.
\end{lem}

We refer to $\Phi_{rep}({\mathcal X})$ as a resolution of ${\mathcal X}$ by representables.  As above, let $R$ be a commutative unital ring.  Using $\Phi_{rep}$, one can define a motive and Suslin homology of any ${\mathcal X} \in \simpnis$.

\begin{defn}
\label{defn:suslinhomology}
Suppose ${\mathcal X}$ is a $k$-space.  The {\em motive of ${\mathcal X}$}, denoted ${\bf M}({\mathcal X})$, is the class of the normalized complex $R_{tr}(\Phi_{rep}({\mathcal X})_\bullet)$ in $\dmeff$.  The $i$-th Suslin homology sheaf of $\mathcal{X}$, denoted $\H_i^S({\mathcal X})$, is defined as $H_i(L_{\aone}R_{tr}(\Phi_{rep}({\mathcal X})_\bullet))$.
\end{defn}

\begin{rem}
In fact, this construction can be extended to a functor $\ho{k} \to \dmeff$; see, e.g., \cite{WeibelRoadmap} for more details regarding this construction.
\end{rem}

While this is not the usual definition of Suslin homology, it coincides with that one in case $k$ is assumed perfect via the following important foundational result.

\begin{thm}[{\cite[Corollary 14.9]{MVW}}]
\label{thm:suslincomplexisaonelocal}
If $k$ is a perfect field, and $X$ is a smooth $k$-scheme, the complex $C_*^S(R_{tr}(X))$ is $\aone$-local.
\end{thm}

\begin{rem}
Suppose $k$ is a field having characteristic $p$.  If $k$ is not perfect, it is not known whether $C_*^S(R_{tr}(X))$ is $\aone$-local.  On the other hand, unpublished work of Suslin establishes that so long as $p$ is invertible in $R$, then Voevodsky's theorem that homotopy invariant presheaves of $R$-modules with transfers have homotopy invariant cohomology \cite[Theorem 13.8]{MVW} still holds.  Using Suslin's result, if $p$ is invertible in $R$, then one can show that $C_*^S(R_{tr}(X))$ is $\aone$-local.  Thus, in this situation, the definition of Suslin homology given above agrees with the usual definition of Suslin homology.
\end{rem}

\subsubsection*{Comparing homology sheaves}
For any $\mathcal{X} \in \simpnis$, recall that $Sing_*^{\aone}({\mathcal X})$ is defined to be the diagonal of the bisimplicial sheaf $(i,j) \mapsto \underline{\hom}(\Delta^{i}_k,\mathcal{X}_j)$, where $\Delta^i_k$ is the algebraic $i$-simplex and we write $\underline{\hom}$ for the internal hom in the category of Nisnevich sheaves of sets.  By construction, there is a canonical morphism $\mathcal{X} \to Sing_*^{\aone}({\mathcal X})$ that is an $\aone$-weak equivalence (see \cite[p. 88]{MV}).  In particular, for any smooth scheme $X$ the morphism $X \to Sing_*^{\aone}(X)$ induces a morphism $C_*(X,R) \to C_*(Sing_*^{\aone}(X),R)$ that becomes an isomorphism after $\aone$-localization.

For any smooth scheme $Y$, there is a canonical monomorphism of sheaves $R(Y) \to R_{tr}(Y)$ (send a morphism $U \to Y$ to the correspondence defined by its graph), and this construction induces a map $N(R(Sing_*^{\aone}(X))) \to C_*^S(R_{tr}(X))$.  Combining this morphism with the discussion of the previous paragraph and applying the $\aone$-localization functor we get morphisms
\[
C_*^{\aone}(X,R) \longrightarrow C_*^{\aone}(Sing_*^{\aone}(X),R) \longrightarrow L_{\aone}(C_*^S(R_{tr}(X)).
\]
We thus obtain a comparison morphism from $\aone$-homology to Suslin homology, and we summarize this construction as follows.

\begin{cor}
\label{cor:aonetosuslin}
For any smooth scheme $X$, there are comparison maps
\[
\H_i^{\aone}(X,R) \longrightarrow \H_i^S(X,R).
\]
induced by the morphism $C_*^{\aone}(X,R) \to L_{\aone}(C_*^S(R_{tr}(X)))$.
\end{cor}

\begin{rem}
With a little more work, one can extend this comparison map to a natural transformation of functors on spaces.
\end{rem}

The following examples show this morphism is {\em not} an isomorphism in general.

\begin{ex}
\label{ex:aonehomologyofgm}
If $(X,x)$ is a pointed smooth scheme, the zeroth Suslin homology sheaf splits $\H_0^S(X) \isomt \Z \oplus \tilde{\H}_0^{S}(X)$.  The morphism of Corollary \ref{cor:aonetosuslin} is compatible with this spitting and there is an induced morphism $\tilde{\H}_0^{\aone}(X) \to \tilde{\H}_0^S(X)$.  In general, this morphism is not an isomorphism, e.g., for $X = \gm$.  Indeed the sheaf of groups $\gm$ is strictly $\aone$-invariant with transfers given by the usual norm map on units.  Lemma \ref{lem:h0universal} shows that the identity map $\gm \to \gm$ induces a homomorphism of strictly $\aone$-invariant sheaves (with transfers) $\tilde{\H}_0^{\aone}(\gm) \to \gm$ (resp. $\tilde{\H}_0^S(\gm) \to \gm$).  Theorem 3.1 of \cite{SuslinVoevodsky} shows the map $\tilde{\H}_0^S(\gm) \to \gm$ is an isomorphism.  More generally, $\tilde{\H}_0^{S}(\gm^{\wedge n})$ is closely related to Milnor K-theory.  On the other hand, $\tilde{\H}_0^{\aone}(\gm^{\wedge n})$ has been computed by Morel (combine \cite[Theorem 2.37]{MField} and \cite[Theorem 4.46]{MField}): when $n \geq 1$, this sheaf is a mixture of Milnor K-theory and Witt groups, which one refers to as Milnor-Witt K-theory.  Similar examples can be constructed from any smooth proper curve.
\end{ex}

\begin{rem}
There is also a ``stabilized" version $\H_0^{s\aone}(X)$ of $\H_0^{\aone}(X)$ where one ``inverts $\gm$."  If $S^0_s = \Spec k_+$, then $\tilde{\H}_0^{s\aone}(S^0_s) = \H_0^{s\aone}(\Spec k)$ essentially by definition.  One can also work directly with the stable $\aone$-homotopy category to prove $\H_0^{s\aone}(\Spec k)$ coincides with the $0$-th stable $\aone$-homotopy sheaf of the motivic sphere spectrum (see \cite[p. 7]{MStable}).  This sheaf has been identified as $\K^{MW}_0$ by \cite[Corollary 6.4.1]{Morelpi0}.

For any finitely generated separable extension $L/k$ there is an isomorphism $\K^{MW}_0(L) \isomt GW(L)$ (see \cite[Remark 6.1.6b]{MStable} or \cite[Lemmas 2.9-2.10]{MField}), where $GW(L)$ denotes the Grothendieck-Witt group of isomorphism classes of non-degenerate symmetric bilinear forms.  On the other hand it follows from the definitions given above that $\H_0^{S}(\Spec k) = \Z$.  These computations suggest that $D_{\aone}(k)$, or perhaps its stabilized version, provides a version of Voevodsky's triangulated category of motives incorporating data from the theory of quadratic forms.  This point of view is further developed in \cite{AH2}.
\end{rem}

\subsubsection*{Connectivity and the $t$-structure}
We now recall some basic facts regarding the structure of $\aone$ (or Suslin) homology sheaves, all due to Morel.  Recall that a complex of sheaves $A_*$ is called {\em $(-1)$-connected (or positive)} if its homology sheaves $H_i(A_*)$ vanish for $i < 0$.  Each of the following results was proven in the context of stable $\aone$-homotopy theory by Morel in \cite{MStable}.  However, the proofs he gives apply just as well (as he observes) to the setting of derived categories that we consider.  For this reason, we give references to the corresponding statements in stable $\aone$-homotopy theory.

\begin{thm}[{\cite[Theorem 6.1.8]{MStable}}]
\label{thm:stableaoneconnectivity}
If $A_*$ is a $(-1)$-connected complex of $R$-modules, then its $\aone$-localization $L_{\aone}(A_*)$ is also $(-1)$-connected.
\end{thm}

\begin{thm}[{\cite[Theorem 6.2.7]{MStable}}]
\label{thm:strictaoneinvarianceofhomology}
If ${\mathcal X}$ is a $k$-space, then for every integer $i$ the sheaves $\H_i^{\aone}({\mathcal X},R)$ and $\H_i^{S}({\mathcal X},R)$ are always strictly $\aone$-invariant, and these sheaves are trivial if $i < 0$.
\end{thm}

\begin{proof}
By construction, the complexes $C_*^{\aone}({\mathcal X},R)$ are $\aone$-localizations of the complexes $C_*({\mathcal X},R)$.  The latter complexes are $(-1)$-connected by definition.  The result follows immediately from Theorem \ref{thm:stableaoneconnectivity}.  The statement for Suslin homology is proven in an identical fashion.
\end{proof}

Before we state the next result, let us recall a variant of \cite[D\'efinition 1.3.1]{BBD}.

\begin{defn}[Homological $t$-structure]
\label{defn:homologicaltstructure}
Let ${\mathcal T}$ be a triangulated category.  A {\em homological $t$-structure} on ${\mathcal T}$ consists of a pair of strictly full subcategories ${\mathcal T}_{\leq 0} \subset {\mathcal T}$ and ${\mathcal T}_{\geq 0} \subset {\mathcal T}$ such that, setting ${\mathcal T}_{\leq n}:= {\mathcal T}_{\leq 0}[n]$ and ${\mathcal T}_{\geq n} := {\mathcal T}_{\geq 0}[n]$, the following properties hold:
\begin{itemize}
\item[i)] for any $C \in {\mathcal T}_{\geq 0}$ and any $D \in {\mathcal T}_{\leq -1}$, one has $\hom_{{\mathcal T}}(C,D) = 0$;
\item[ii)] there are inclusions ${\mathcal T}_{\geq 1} \subset {\mathcal T}_{\geq 0}$, and ${\mathcal T}_{\leq 0} \subset {\mathcal T}_{\leq 1}$;
\item[iii)] for any object $X \in {\mathcal T}$ there is a distinguished triangle
\[
A \longrightarrow X \longrightarrow B \longrightarrow A[1]
\]
with $A \in {\mathcal T}_{\geq 0}$ and $B \in {\mathcal T}_{\leq 1}$.
\end{itemize}
\end{defn}

If $({\mathcal T},{\mathcal T}_{\leq 0},{\mathcal T}_{\geq 0})$ is a homological $t$-structure on ${\mathcal T}$, then by \cite[Proposition 1.3.3]{BBD} there are truncation functors $\tau_{\geq n}: {\mathcal T} \to {\mathcal T}_{\geq n}$, and $\tau_{\leq n}: {\mathcal T} \to {\mathcal T}_{\leq n}$ adjoint to the corresponding inclusion functors.  Using homological conventions as above, the proof of {\em loc. cit.} gives the following result.

\begin{prop}
\label{prop:homologicaltstructure}
Suppose $({\mathcal T},{\mathcal T}_{\leq 0},{\mathcal T}_{\geq 0})$ is a $t$-category.  If $X$ is any object in ${\mathcal T}$, there exists a unique morphism $d \in \hom^1(\tau_{\leq 1}X,\tau_{\geq 0}X)$ such that the triangle
\[
\tau_{\geq 0} X \longrightarrow X \longrightarrow \tau_{\leq 1} X \stackrel{d}{\longrightarrow} \tau_{\geq 0}X[1]
\]
is distinguished.
\end{prop}

A complex $A_*$ is called {\em negative} if $H_i(A_*) = 0$ for $i > 0$ and {\em positive} if $H_i(A_*) = 0$ for $i < 0$.  We write $D_{\aone}(k)_{\leq 0}$ for the full subcategory of $D_{\aone}(k)$ consisting of $\aone$-local negative complexes, and $D_{\aone}(k)_{\geq 0}$ for the full subcategory consisting of $\aone$-local positive complexes.

\begin{prop}[{\cite[Lemma 6.2.11]{MStable}}]
\label{prop:aonetstructure}
The triple $(D_{\aone}(k),D_{\aone}(k)_{\leq 0},D_{\aone}(k)_{\geq 0})$ is a homological $t$-structure on $D_{\aone}(k)$.
\end{prop}

If we write $\Ab^{\aone}_k$ for the category of strictly $\aone$-invariant sheaves.  The category of strictly $\aone$-invariant sheaves of groups can be identified as the heart of this $t$-structure. By \cite[Th\'eor\`eme 1.3.6]{BBD}, we get the following result.

\begin{cor}
\label{cor:strictlyaoneinvariantabeliancategory}
The category $\Ab^{\aone}_k$ is abelian.
\end{cor}

If $k$ is perfect, Voevodsky showed that $\dmeff$ admits a $t$-structure defined in a manner identical to Proposition \ref{prop:aonetstructure}.  The heart of the resulting $t$-structure is precisely the category of strictly $\aone$-invariant Nisnevich sheaves {\em with transfers}; see \cite[\S 4.3]{Deglise1} for a discussion.  We write $\Ab^{\aone}_{tr,k}$ for the abelian category of strictly $\aone$-invariant sheaves with transfers.  We can define $\H_0^{S}$ as a functor from $\dmeff$ to $\Ab^{\aone}_{tr,k}$ (see {\em loc. cit.} Formula 4.12a).

\subsubsection*{Gersten resolutions}
Suppose $A$ is an $\aone$-local complex.  The axiomatic approach of \cite{CTHK} provides general machinery for producing a Gersten resolution associated with the Nisnevich hypercohomology of $A$ (see {\em ibid.} \S 7).  Let $(X,Z)$ be a pair where $X$ is a smooth scheme and $Z \subset X$ is a closed subscheme.  The functor
\[
(X,Z) \longmapsto {\mathbb H}^*_{Z}(X_{Nis},A)
\]
(i.e., Nisnevich hypercohomology with supports on $Z$) defines a cohomology theory with supports in the sense of \cite[Definition 5.1.1]{CTHK}.  Moreover, this theory satisfies Nisnevich excision (Axiom {\bf COH1} of \cite[p. 55]{CTHK}) and $\aone$-homotopy invariance (Axiom {\bf COH3} of {\em ibid} p. 58).  Let ${\mathbb H}^n_{Zar}(A)$ denote the Zariski sheaf associated with the presheaf $U \mapsto {\mathbb H}^n_{Nis}(U,A)$; this apparent abuse of terminology will be justified momentarily.

\begin{prop}[{\cite[Corollary 5.1.11]{CTHK}}]
\label{prop:gerstenresolution}
Suppose $k$ is an infinite field.  For any smooth $k$-scheme $X$, and any $\aone$-local complex $A$, the complex
\[
{\mathbb H}^n_{Zar}(A)|_X \longrightarrow \coprod_{x \in X^{(0)}} i_{x_*}{\mathbb H}^n_x(X,A) \longrightarrow \cdots \longrightarrow \coprod_{x \in X^{(p)}} i_{x_*}{\mathbb H}^{p+n}_x(X,A) \longrightarrow \cdots
\]
is a flasque resolution.
\end{prop}

\begin{proof}
We just observe that \cite{CTHK} Proposition 5.3.2a, Axiom {\bf COH2} is implied by Axiom {\bf COH3}).
\end{proof}

We write ${\mathbb H}^n_{Nis}(A)$ for the {\em Nisnevich} sheaf associated with the presheaf $U \mapsto {\mathbb H}^n_{Nis}(U,A)$.  We use the following fundamental comparison result.

\begin{thm}[{\cite[Theorem 8.3.1]{CTHK}}]
\label{thm:comparisonZariskiNisnevich}
Suppose $k$ is an infinite field.  For any smooth $k$-scheme $X$, and any $\aone$-local complex $A$, the canonical maps
\[
H^i_{Zar}(X,{\mathbb H}^n_{Zar}(A)) \longrightarrow H^i_{Nis}(X,{\mathbb H}^n_{Nis}(A))
\]
are isomorphisms.
\end{thm}

\section{Birational geometry and strictly $\aone$-invariant sheaves}
\label{s:birationalproperties}
In this section, we study the relationship between the zeroth $\aone$-homology or Suslin homology sheaf (recall Definitions \ref{defn:aonehomology} and \ref{defn:suslinhomology}) and the zeroth $\aone$-homotopy sheaf (recall Definition \ref{defn:aoneconnected}).  The principal results of this section imply Proposition \ref{propintro:universalproperty} and Theorem \ref{thmintro:birationalinvariance} of the introduction.  We prove in Lemma \ref{lem:h0universal} that the zeroth $\aone$-homology (resp. Suslin homology) sheaf of a smooth scheme $X$ is initial among strictly $\aone$-invariant sheaves (with transfers) admitting a morphism from $X$. In Proposition \ref{prop:dependsonpi0} we establish that the zeroth $\aone$-homology (resp. Suslin homology) sheaf of $X$ is the free strictly $\aone$-invariant sheaf (with transfers) generated by $\pi_0^{\aone}(X)$, and in Theorem \ref{thm:birationalclass} that it is a stable birational invariant for smooth proper schemes over infinite fields.

\subsubsection*{A factorization lemma}
The functor sending a Nisnevich sheaf $\F$ to the corresponding constant simplicial sheaf (all face and degeneracy maps are the identity) is fully-faithful.  The full-subcategory of $\simpnis$ consisting of constant simplicial sheaves will be referred to as the subcategory of spaces of simplicial dimension $0$ (see \cite[p. 47]{MV}).  Spaces of simplicial dimension $0$ are automatically simplicially fibrant (see \cite[\S 2 Remark 1.14]{MV}); in particular, since smooth schemes have simplicial dimension $0$, they are simplicially fibrant.  If ${\mathcal X}$ is any space, then the unstable $\aone$-$0$-connectivity theorem \cite[\S 2 Corollary 3.22]{MV} gives an epimorphism ${\mathcal X} \to \pi_0^{\aone}({\mathcal X})$.  We begin by stating a result that will be used repeatedly in the sequel.

\begin{lem}
\label{lem:factorization}
Suppose ${\mathcal X}$ and ${\mathcal Y}$ are $k$-spaces of simplicial dimension $0$, and ${\mathcal Y}$ is $\aone$-local.  The canonical epimorphism ${\mathcal X} \to \pi_0^{\aone}({\mathcal X})$ induces a bijection
\begin{equation}
\label{eqn:factorization}
\hom_{\Shv_{Nis}(\Sm_k)}(\pi_0^{\aone}({\mathcal X}),{\mathcal Y}) \longrightarrow \hom_{\Shv_{Nis}(\Sm_k)}({\mathcal X},{\mathcal Y})
\end{equation}
functorial in both inputs.
\end{lem}

\begin{proof}
Since ${\mathcal X}$ and ${\mathcal Y}$ have simplicial dimension $0$, the canonical map
\[
\hom_{\Shv_{Nis}(\Sm_k)}({\mathcal X},{\mathcal Y}) \longrightarrow [{\mathcal X},{\mathcal Y}]_s
\]
is a bijection, as we observed just before the statement of the lemma.  Since ${\mathcal Y}$ is $\aone$-local, we also have identifications $[{\mathcal X},{\mathcal Y}]_{\aone} = [{\mathcal X},{\mathcal Y}]_s$.

Since smooth schemes have simplicial dimension $0$, for any smooth scheme $U$, we have $[U,{\mathcal Y}]_{\aone} = \hom_{\Shv_{Nis}(\Sm_k)}(U,{\mathcal Y})$.  Sheafifying for the Nisnevich topology, we deduce that $\pi_0^{\aone}({\mathcal Y}) = {\mathcal Y}$.  To finish, observe that any morphism ${\mathcal X} \to \pi_0^{\aone}({\mathcal Y})$ factors uniquely through $\pi_0^{\aone}({\mathcal X})$ by the definition of $\pi_0^{\aone}(\cdot)$.
\end{proof}

\begin{cor}
\label{cor:aonelocaldetection}
Suppose $X$ is an $\aone$-connected smooth $k$-scheme.  If ${\mathcal Y}$ is an $\aone$-local space of simplicial dimension $0$, then the map ${\mathcal Y}(k) \to {\mathcal Y}(X)$ induced by the structure map is a bijection.  In particular, if $M$ is a strictly $\aone$-invariant sheaf, and $X$ is an $\aone$-connected smooth $k$-scheme, the canonical map $M(k) \to M(X)$ is a bijection.
\end{cor}

\begin{lem}
\label{lem:h0universal}
If $M$ (resp. $M'$) is a strictly $\aone$-invariant sheaf of $R$-modules (with transfers), then for any space ${\mathcal X}$ there are bijections
\[
\begin{split}
H^0_{Nis}({\mathcal X},M) &\isomto \hom_{\Ab^{\aone}_k}(\H_0^{\aone}({\mathcal X},R),M) \\
H^0_{Nis}({\mathcal X},M') &\isomto \hom_{\Ab^{\aone}_k}(\H_0^{S}({\mathcal X},R),M')
\end{split}
\]
functorial in both ${\mathcal X}$ and $M$ (resp. $M'$).
\end{lem}

\begin{proof}
As above, since $M$ is $\aone$-local we have identifications
\[
\hom_{\Shv_{Nis}(\Sm_k)}({\mathcal X},M) \isomto [{\mathcal X},M]_{\aone} = [{\mathcal X},K(M,0)]_{\aone}.
\]
The adjunction between the $\aone$-homotopy and $\aone$-derived categories allows one to identify the last abelian group with $\hom_{D_{\aone}(k)}(C_*^{\aone}({\mathcal X},R),M)$.

Now, by Proposition \ref{prop:aonetstructure}, we know that $D_{\aone}(k)$ admits a homological $t$-structure.  By the stable $\aone$-connectivity theorem, we know that $C_*^{\aone}({\mathcal X},R) \in D_{\aone}(k)_{\geq 0}$.  The result follows from a general fact about homological $t$-structures (see Definition \ref{defn:homologicaltstructure}). Let $\mathcal{T}$ be a triangulated category with a homological $t$-structure $(\mathcal{T}_{\geq 0},\mathcal{T}_{\leq 0})$ and heart $\mathcal{A}$.  Let $M \in {\mathcal A}$, and write $M[0]$ for $M$ viewed as an object of ${\mathcal T}$ situated in degree $0$.  For any object $C \in {\mathcal T}_{\geq 0}$, Proposition \ref{prop:homologicaltstructure} and a shifting argument give rise to the distinguished triangle
\[
\tau_{\geq -1}C \longrightarrow C \longrightarrow \tau_{\leq 0}C \longrightarrow \tau_{\geq -1}C[1].
\]
Applying the functor $\hom_{{\mathcal T}}(\cdot,M[0])$ to this distinguished triangle, we get a map
\[
\hom_{{\mathcal T}}(\tau_{\leq 0}C,M[0]) \longrightarrow \hom_{{\mathcal T}}(C,M[0]),
\]
and this map is an isomorphism directly from Definition \ref{defn:homologicaltstructure}(i) and (ii), which show that $\hom_{{\mathcal T}}(\tau_{\geq -1}C,M[0])$ and $\hom_{{\mathcal T}}(\tau_{\geq -1}C[1],M[0])$ vanish.  Since $\tau_{\geq 0}C = C$, and $H_0(C) = \tau_{\leq 0}\tau_{\geq 0}C$ (see \cite[Th\'eor\`eme 1.3.6]{BBD}) we get an isomorphism
\[
\hom_{{\mathcal A}}(H_0(C),M) \isomto \hom_{\mathcal{T}}(C,M[0]).
\]

The proof for Suslin homology sheaves is similar; one uses \cite[Exercise 13.6]{MVW} to identify $\hom_{\Shv_{Nis}(\Sm_k)}({\mathcal X},M')$ with $\hom_{\dmeff}({\mathbf M}({\mathcal X}),M)$ (since $M'$ is $\aone$-local) together with an identical truncation argument.
\end{proof}

\begin{cor}
\label{cor:factorization}
If $X$ is a smooth $k$-scheme, and $M$ is a strictly $\aone$-invariant sheaf of $R$-modules with transfers, then any morphism $\varphi: \H_0^{\aone}(X,R) \to M$ factors as a composite
\[
\H_0^{\aone}(X,R) \longrightarrow \H_0^{S}(X,R) \longrightarrow M,
\]
where the first map is the morphism of \textup{Corollary \ref{cor:aonetosuslin}}.
\end{cor}

\begin{proof}
The morphism of Corollary \ref{cor:aonetosuslin} induces for any strictly $\aone$-invariant sheaf of $R$-modules with transfers a morphism
\[
\hom_{\Ab^{\aone}_k}(\H_0^{S}(X,R),M) \longrightarrow \hom_{\Ab^{\aone}_k}(\H_0^{\aone}(X,R),M).
\]
Lemma \ref{lem:h0universal} implies that this morphism is a bijection.
\end{proof}

\subsubsection*{Dependence on the sheaf of $\aone$-connected components}
If $M$ is a (sufficiently nice) topological space, the Mayer-Vietoris sequence shows that the ordinary singular homology group $H_0(M,R)$ is the free $R$-module generated by the connected components of $M$.  On the other hand, Lemma \ref{lem:h0universal} can be interpreted as saying that $\H_0^{\aone}({\mathcal X},R)$ (resp. $\H_0^{S}({\mathcal X},R)$) is the free strictly $\aone$-invariant sheaf (with transfers) on the sheaf ${\mathcal X}$.  We now show that $\H_0^{\aone}({\mathcal X},R)$ (resp. $\H_0^{S}({\mathcal X},R)$) is the free strictly $\aone$-invariant sheaf (with transfers) on the sheaf of $\aone$-connected components of ${\mathcal X}$.

\begin{prop}
\label{prop:dependsonpi0}
For any space ${\mathcal X}$, and any commutative unital ring $R$, the maps
\[
\begin{split}
\H_0^{\aone}({\mathcal X},R) &\longrightarrow \H_0^{\aone}(\pi_0^{\aone}({\mathcal X}),R) \text{ and }\\
\H_0^{S}({\mathcal X},R) &\longrightarrow \H_0^{S}(\pi_0^{\aone}({\mathcal X}),R),
\end{split}
\]
induced by the canonical epimorphism ${\mathcal X} \to \pi_0^{\aone}({\mathcal X})$ are isomorphisms.
\end{prop}

\begin{proof}
We prove the first statement; the second statement is proven in an essentially identical manner.  Assume first that ${\mathcal X}$ has simplicial dimension $0$. Let $M$ be an arbitrary strictly $\aone$-invariant sheaf (with transfers for the second statement).  We have a commutative diagram
\[
\xymatrix{
\hom_{\Shv_{Nis}(\Sm_k)}(\pi_0^{\aone}({\mathcal X}),M) \ar[r]\ar[d] & \hom_{\Ab^{\aone}_k}(\H_0^{\aone}(\pi_0^{\aone}({\mathcal X})),M) \ar[d] \\
\hom_{\Shv_{Nis}(\Sm_k)}({\mathcal X},M) \ar[r] & \hom_{\Ab^{\aone}_k}(\H_0^{\aone}({\mathcal X}),M)
}.
\]
By Lemma \ref{lem:h0universal} the horizontal maps are isomorphisms, and by Lemma \ref{lem:factorization} the left vertical map is a bijection.  Indeed, all these bijections are functorial in both variables.  It follows that the right vertical map is a bijection functorially in both variables as well.  The result then follows from the Yoneda lemma.

To treat the general case, it suffices to observe that by \cite[\S 2 Proposition 3.14]{MV} every space ${\mathcal X}$ is $\aone$-weakly equivalent to a space of simplicial dimension $0$.
\end{proof}

\begin{rem}[A non-abelian variant]
One may also prove a non-abelian version of Proposition \ref{prop:dependsonpi0}. Because one needs to keep track of base points, this version seems not as widely applicable.  Recall that if $({\mathcal S},s)$ is a pointed sheaf of sets, we can consider $F_{\aone}({\mathcal S}) := \pi_1^{\aone}(\Sigma^1_s {\mathcal S})$.  Results of \cite{MField} show that this sheaf is strongly $\aone$-invariant (see Remark \ref{rem:stronginvariance}).  As above, one can show that the canonical map $F_{\aone}({\mathcal S}) \to F_{\aone}(\pi_0^{\aone}({\mathcal S}))$ is an isomorphism.  This result is compatible with the previous results via the $\aone$-Hurewicz theorem (also proven by Morel).  The sheaves $F_{\aone}(\pi_0^{\aone}({\mathcal S}))$ contain ``non-abelian" information, e.g., related to finite covers with non-abelian fundamental group.
\end{rem}

\begin{lem}[{\cite[Lemma 6.4.4]{MStable}}]
\label{lem:excision}
Suppose $M$ is a strictly $\aone$-invariant sheaf of groups.  If $X$ is a smooth $k$-scheme, and $U \subset X$ is an open subscheme whose complement has codimension $\geq d$ in $X$, then the restriction map
\[
H^i_{Nis}(X,M) \to H^i_{Nis}(U,M)
\]
is a monomorphism if $i \leq d-1$ and a bijection if $i \leq d-2$.
\end{lem}

\begin{prop}
\label{prop:epimorphism}
Suppose $X$ is a smooth $k$-scheme and $U \subset X$ is an open subscheme of $X$.  Assume the complement of $U$ in $X$ has codimension $\geq d$, for some integer $d > 0$.  For any commutative unital ring $R$, the canonical maps
\[
\begin{split}
\H_0^{\aone}(U,R) &\longrightarrow \H_0^{\aone}(X,R), \text{ and } \\
\H_0^{S}(U,R) &\longrightarrow \H_0^{S}(X,R)
\end{split}
\]
are epimorphisms if $d = 1$, and isomorphisms if $d \geq 2$.
\end{prop}

\begin{proof}
For the first statement, if $M$ is an arbitrary strictly $\aone$-invariant sheaf of $R$-modules, we have functorial bijections $\hom_{\Ab^{\aone}_k}(\H_0^{\aone}(X,R),M) \isomt M(X)$ by Lemma \ref{lem:h0universal}.  Likewise, if $M$ is an arbitrary strictly $\aone$-invariant sheaf of $R$-modules with transfers, we have functorial bijections $\hom_{\Ab^{\aone}_{tr,k}}(\H_0^{S}(X,R),M) \isomt M(X)$ by Lemma \ref{lem:h0universal}. Thus, the result follows immediately from Lemma \ref{lem:excision} and the Yoneda lemma.
\end{proof}

\subsubsection*{Stable birational equivalence}
Recall that two smooth proper $k$-varieties $X$ and $Y$ are stably $k$-birationally equivalent if $X \times {\mathbb P}^n$ is $k$-birationally equivalent to $Y \times {\mathbb P}^m$ for integers $m,n \geq 0$.  In particular, if $X$ is stably $k$-birationally equivalent to projective space, then we say that $X$ is stably $k$-rational.

\begin{thm}
\label{thm:birationalclass}
Suppose $k$ is an infinite field, and $R$ is a commutative unital ring.  If $X$ and $X'$ are stably $k$-birationally equivalent smooth proper varieties then $\H_0^{\aone}(X,R) \cong \H_0^{\aone}(X',R)$ and $\H_0^{S}(X,R) \cong \H_0^{S}(X',R)$.
\end{thm}

\begin{proof}
Consider the composite map $Y \times {\mathbb A}^n \hookrightarrow Y \times {\mathbb P}^{n} \longrightarrow Y$.  Since $\H_0^{\aone}(Y)$ is $\aone$-homotopy invariant, it follows that the composite map $\H_0^{\aone}(Y \times {\mathbb A}^n) \longrightarrow \H_0^{\aone}(Y)$ is an isomorphism.  On the other hand, the map $\H_0^{\aone}(Y \times {\mathbb A}^n) \longrightarrow \H_0^{\aone}(Y \times {\mathbb P}^n)$ is an epimorphism by Proposition \ref{prop:epimorphism}.  A diagram chase shows that projection map must then also be an isomorphism.  The same argument works for Suslin homology.

If $k$ has characteristic $0$, we may finish the proof by means of a straightforward geometric argument using resolution of singularities. Indeed, given any $k$-birational morphism $X \to Y$, there is a commutative diagram of $k$-birational morphisms of the form
\[
\xymatrix{
X' \ar[r]\ar[d] & Y' \ar[d]\ar[dl] \\
X \ar[r] & Y
}
\]
where all the vertical maps are composites of a finite number of blow-ups with smooth centers.  We claim that it suffices to show that the morphism on zeroth $\aone$-homology sheaves induced by a blow-up with smooth center is an isomorphism.  If that is the case, since the composite map $\H_0^{\aone}(X') \to \H_0^{\aone}(Y') \to \H_0^{\aone}(X)$ is an isomorphism, we realize $\H_0^{\aone}(X)$ as a summand of $\H_0^{\aone}(Y)$ and vice versa (by reversing the roles of $X$ and $Y$).

Let us check the result for $f: X' \to X$, where $f$ is a blow-up at a codimension $\geq 2$ smooth subscheme $Z \subset X$.  The induced map $X' \setminus f^{-1}(Z) \to X \setminus Z$ is an isomorphism.  The morphism $X' \setminus f^{-1}(Z) \to X'$ is an open immersion with complement having codimension $1$, and the morphism $X \setminus Z \to X$ is an open immersion with complement having codimension $\geq 2$.  By the previous proposition, the map $\H_0^{\aone}(X' \setminus f^{-1}(Z),R) \to \H_0^{\aone}(X',R)$ is an epimorphism, the map $\H_0^{\aone}(X' \setminus f^{-1}(Z),R) \to \H_0^{\aone}(X \setminus Z,R)$ is an isomorphism, and the map $\H_0^{\aone}(X \setminus Z,R) \to \H_0^{\aone}(X,R)$ is an isomorphism.  Composing the second and third of these isomorphisms, we see that the morphism $\H_0^{\aone}(X' \setminus f^{-1}(Z),R) \to \H_0^{\aone}(X,R)$ is an isomorphism.  Consequently, the morphism $\H_0^{\aone}(X',R) \to \H_0^{\aone}(X,R)$ is an isomorphism as well.  The case of the zeroth Suslin homology sheaf is identical.

If $k$ is just infinite, we argue as follows.  If $M$ is an arbitrary strictly $\aone$-invariant sheaf, then $M$ is $\aone$-local and so admits a Gersten resolution by Proposition \ref{prop:gerstenresolution}.  By \cite[Theorem 8.5.1]{CTHK} it follows that $M(X)$ is a birational invariant of smooth proper varieties (this explanation is expanded slightly in Lemma \ref{lem:unramifiedinvariants}).  Since $M$ was arbitrary, it follows from \ref{lem:h0universal} that the same statement holds for the zeroth $\aone$-homology sheaf.  An analogous argument works for the zeroth Suslin homology sheaf.
\end{proof}

\begin{rem}
As discussed in \cite[\S 2]{AM}, we know that $\aone$-connectedness is a birational invariant for fields having characteristic $0$.  However, we do not know whether the sheaf $\pi_0^{\aone}(X)$ is itself a (stable) birational invariant. The above result shows that after abelianization this is the case.  Note also that the first proof above implies that if $X$ and $X'$ are two schemes, not necessarily proper, that can be linked by a chain of blow-ups at smooth schemes, then $\H_0^{\aone}(X) \cong \H_0^{\aone}(X')$.  In fact, it is not at the moment known whether $\pi_0^{\aone}(X)$ is unchanged by blow-ups along smooth schemes!
\end{rem}

\begin{rem}
\label{rem:birationalinvarianceofchow}
In Example \ref{ex:unramifiedchow} we will see that if $k$ is a perfect field and $X$ is a smooth proper $k$-variety, then for any separable finitely generated extension $L/k$, one can identify $\H_0^S(X)(L)$ with $CH_0(X_L)$, functorially in $L$.  However, birational invariance for the Chow group of $0$-cycles is a much older result.  Indeed, if $k$ has characteristic $0$ then \cite[Proposition 6.3]{CTC} establishes $k$-birational invariance of $CH_0(X)$, and Fulton \cite[Example 16.1.11]{Fulton} generalizes this to arbitrary characteristic.
\end{rem}

\section{An unramified characterization of $\aone$-connectedness}
\label{s:unramifiedelements}
In this section, we recall aspects of a ``field theoretic" or ``unramified" approach to strictly $\aone$-invariant sheaves pioneered by Morel \cite{MMilnor,MStable,MField} following foundational work of Rost \cite{RostChow}.  If $M$ is a strictly $\aone$-invariant sheaf (see Definition \ref{defn:strictaoneinvariance}), we now explain how to identify sections of $M$ over a smooth scheme $X$ in terms of the function field $k(X)$ of $X$ and codimension $1$ geometry of $X$, i.e., geometric discrete valuations on $k(X)$.  We then provide a number of examples of strictly $\aone$-invariant sheaves.  Combining this result with the discussion of \S \ref{s:birationalproperties} (specifically Lemma \ref{lem:h0universal}) Lemma \ref{lem:unramifiedinvariants} explains the sense in which $\H_0^{\aone}(X)$ is a ``universal unramified invariant" as mentioned in the introduction.  Theorem \ref{thm:characterization} and the subsequent corollary give the the unramified characterization of $\aone$-connectedness stated in the introduction.  By means of an example, we show that Suslin homology is not sufficiently refined to detect $\aone$-connectedness, or more loosely that $\aone$-connectedness cannot be characterized solely by means of ``unramified invariants with transfers;"  see Proposition \ref{prop:rationallyconnected}, Example \ref{ex:artinmumford} and Example \ref{ex:parimala} for more details.

\subsubsection*{Unramified elements}
Fix a field $k$, and suppose $M$ is a strictly $\aone$-invariant sheaf (on $\Sm_k$). Suppose $S$ is an essentially smooth $k$-scheme, i.e., a filtering inverse limit of smooth schemes with smooth affine transition morphisms.  If we write $S = \lim X_{\alpha}$, we can define $M(S) = \colim M(X_{\alpha})$.  One can check that this colimit is independent of the choice of filtering inverse system defining $S$.  Thus, we can extend $M$ uniquely to a functor on the category of essentially smooth $k$-schemes.  If $\F_k$ denotes the category of finitely generated extension fields (morphisms are inclusions of fields), then $M$ gives rise to a (covariant) functor on $\F_k$.  Abusing notation, we will denote all these functors by $M$.

By Lemma \ref{lem:excision}, for an open immersion of smooth schemes $U \hookrightarrow X$, the restriction map $M(X) \to M(U)$ is injective.  If $L/k$ is a finitely generated extension of $k$, $\nu$ is a geometric discrete valuation of $L$ with valuation ring $\O_{\nu}$, and $\kappa_{\nu}$ is the associated residue field then we have a morphism $M(\O_{\nu}) \to M(L)$; this morphism is injective by what we've just said.  We now use these observations to define unramified groups associated with any strictly $\aone$-invariant sheaf.

\begin{defn}
Suppose $X$ is an irreducible smooth $k$-scheme.  Given $x \in X^{(1)}$, write $\nu_{x}$ for the corresponding discrete valuation.  For any $x \in X^{(1)}$, the map $M(X) \to M(\O_{\nu_x})$ is injective.  Set
\[
M^{ur}(X) := \bigcap_{x \in X^{(1)}} M(\O_{\nu_x}),
\]
where the intersection is taken in $M(k(X))$.
\end{defn}

\begin{lem}
\label{lem:unramifiedinvariants}
The induced map $M(X) \to M^{ur}(X)$ is an isomorphism.  Thus, if $X$ is a smooth variety, the functor $M \mapsto M^{ur}(X)$ (from the category of strictly $\aone$-invariant sheaves of groups to the category of abelian groups) is representable on $\Ab^{\aone}_k$ by the sheaf $\H_0^{\aone}(X)$.
\end{lem}

\begin{proof}
A class $\alpha \in M^{ur}(X)$ comes from a class $M(k(X))$ lying in the image of $M(\O_{\nu_x})$ as $x$ ranges over the codimension $1$ points of $X$.  If $\alpha$ is in the image of $\O_{\nu}$, then by definition there is an open subscheme $U_{\nu} \subset X$ on which $\alpha$ is defined.  Thus, we can find a collection of open subschemes $U_{i}$ such that $\alpha$ extends to a class on $U_i$ for each $i$.  Using the sheaf property and induction, these classes glue to give a class on the union $U$ of the $U_i$.  By assumption, this union contains all codimension $1$ points of $X$.  By Lemma \ref{lem:excision}, we know that if $U \subset X$ is an open subscheme whose closed complement has codimension $\geq 2$, the restriction map $M(X) \to M(U)$ is an isomorphism.  The second statement follows immediately from the first one via Lemma \ref{lem:h0universal}.
\end{proof}

\begin{cor}
\label{cor:isomorphismsofstrictlyaoneinvariantsheaves}
Given $M,M' \in \Ab^{\aone}_k$, then $f: M \to M'$ is an isomorphism if and only if for every separable, finitely generated extension $L/k$ the morphism $M(L) \to M'(L)$ is an isomorphism.
\end{cor}

\begin{proof}
Since $\Ab^{\aone}_k$ is abelian, it suffices to prove that the strictly $\aone$-invariant sheaves $\ker(f)$ and $\operatorname{coker}(f)$ are trivial.  However, it follows immediately from Lemma \ref{lem:unramifiedinvariants} that a strictly $\aone$-invariant sheaf $A$ is trivial if and only if $A(L)$ is trivial.
\end{proof}

\begin{rem}
If $M$ admits transfers, then the point of view on strictly $\aone$-invariant sheaves (with transfers) discussed above is closely related to Rost's theory of cycle modules \cite{RostChow}.  In fact, all of the examples of strictly $\aone$-invariant sheaves used below can be constructed using either Rost's theory or a modification developed by Morel.  The relationship between strictly $\aone$-invariant sheaves with transfers and Rost's theory of cycle modules has been developed by D\'eglise \cite{Deglise1,Deglise2} (the former category is a localization of the latter).  The counterpart of Lemma \ref{lem:unramifiedinvariants} in the setting of cycle modules is given by a result of Merkurjev \cite[Theorem 2.10]{Merkurjev}.  In fact, the strictly $\aone$-invariant sheaf with transfers associated with Merkurjev's cycle module by D\'eglise's theory (or rather its degree $0$ part) is precisely the $0$-th Suslin homology sheaf as we explain below in Example \ref{ex:unramifiedchow}.
\end{rem}

\begin{rem}
There is a quotient map $\O_{\nu} \to \kappa_{\nu}$, and this induces a morphism $M(\O_{\nu}) \to M(\kappa_{\nu})$.  By choosing local parameters, one can define appropriate notions of residue maps for strictly $\aone$-invariant sheaves, though if $M$ does not admit transfers, these residues depend on the choices made.  This point of view is developed in \cite[\S 1]{MField}, but we will not use this theory below.
\end{rem}

\subsubsection*{Unramified \'etale cohomology and other examples}
Recall by Corollary \ref{cor:aonelocaldetection}, if $M$ is a strictly $\aone$-invariant sheaf, and $X$ is an $\aone$-connected smooth scheme, the pullback map $M(\Spec k) \to M(X)$ is a bijection.  In this section, we give a number of examples of unramified sheaves to show what kind of ``vanishing" statements $\aone$-connectedness entails.

\begin{ex}
\label{ex:unramifiedetalecohomology}
Suppose $k$ is a field, and $n$ is an integer that is not divisible by the characteristic of $k$.  Let ${\mathbf H}^p_{\et}(\mu_n^{\tensor q})$ denote the (Nisnevich) sheaf (on $\Sm_k$) associated with the presheaf $U \mapsto H^p_{\et}(U,\mu_n^{\tensor q})$.  The sheaf ${\mathbf H}^p_{\et}(\mu_n^{\tensor q})$ is strictly $\aone$-invariant.  There are many ways to see this; for example, it follows with a bit of work from homotopy invariance for \'etale cohomology (\cite[Expose XV Lemme 4.2]{SGA43}).  %Since we will not use it directly, we just hit it with a sledgehammer: Corollary \ref{cor:sheafifiedblochkato}.
\end{ex}

\begin{lem}
\label{lem:unramifiedrelation}
Suppose $n$ is an integer that is coprime to the characteristic of $k$.  If $X \in \Sm_k$, then we have
\[
\hom_{\Ab^{\aone}_k}(\H_0^{\aone}(X),{\mathbf H}^p_{\et}(\mu_n^{\tensor q})) = H^0_{Nis}(X,{\mathbf H}^p_{\et}(\mu_n^{\tensor q})).
\]
If furthermore $X$ is proper, then the latter group is precisely the group $H^p_{ur}(k(X)/k,\mu_n^{\tensor q})$.
\end{lem}

\begin{proof}
Since ${\mathbf H}^p_{\et}(\mu_n^{\tensor q})$ is strictly $\aone$-invariant, the equality in the statement follows immediately from Lemma \ref{lem:h0universal}.  By Lemma \ref{lem:unramifiedinvariants}, if $X$ is an irreducible smooth scheme, ${\mathbf H}^p_{\et}(\mu_n^{\tensor q})(X)$ coincides with the subgroup of $H^p_{\et}(k(X),\mu_n^{\tensor q})$ consisting of unramified elements.
%The last statement follows, e.g., from the equivalence of several definitions of unramified (\'etale) cohomology given in \cite[Theorem 4.1.1]{CTPurity}.
\end{proof}

\begin{rem}
For a development of unramified \'etale cohomology see, e.g., \cite[\S 4]{CTPurity}.  These groups admit an alternate description.  If $L/k$ is a finitely generated field extension and $\nu$ is a discrete valuation of $L/k$ with residue field $\kappa$, there are residue maps
\[
\partial_{\nu}: H^p_{\et}(L,\mu_n^{\tensor q}) \longrightarrow H^{p-1}_{\et}(\kappa,\mu_n^{\tensor q-1}).
\]
The subgroup of $H^p_{\et}(L,\mu_n^{\tensor q})$ can be identified with the intersection of the kernels of the residue maps $\partial_{\nu}$ as $\nu$ ranges over the discrete valuations of $L/k$.
\end{rem}

\begin{ex}
\label{ex:unramifiedMilnorktheory}
Let $X$ be a smooth $k$-variety, with function field $k(X)$.  Given a codimension $1$ point, we write $\partial_x$ for the residue map associated with the valuation ring defined by $x$ (see \cite[Lemma 2.1]{Milnor} for the construction of these residue maps).  We define
\[
\K^M_n(X) := \ker(K^M_n(k(X)) \stackrel{\bigoplus_{x \in X^{(1)}}\partial_x}{\longrightarrow} \bigoplus_{x \in X^{(1)}} K^M_{n-1}(\kappa_{\nu})).
\]
We recall one functoriality property of these residue maps immediately subsequent to this example.  If $\varphi: A \to B$ is a ring homomorphism, the induced map $K^M_n(A) \to K^M_n(B)$ is usually denoted by either $\varphi_*$ or $Res_{A/B}$.

One can show that $\K^M_n$ is a strictly $\aone$-invariant sheaf (see \cite[\S 2.2]{MMilnor} for more details), and in fact $\K^M_n$ is a strictly $\aone$-invariant sheaf with transfers.  For any integer $m$, multiplication by $m$ extends to a morphism of sheaves $\K^M_n \stackrel{\times m}{\longrightarrow} \K^M_n$.  The category of strictly $\aone$-invariant sheaves is abelian (see Corollary \ref{cor:strictlyaoneinvariantabeliancategory}), and the cokernel of this morphism of sheaves, which is necessarily strictly $\aone$-invariant, is denoted $\K^M_n/m$.  By construction, we have identifications $\K^M_n(L) = K^M_n(L)$ and $\K^M_n/m(L) = K^M_n(L)/m$ for any finitely generated extension $L/k$.
\end{ex}

\begin{ex}
\label{ex:unramifiedchow}
Let $k$ be a perfect field and let $X$ be a smooth proper $k$-variety.  Consider the functor assigning to a finitely generated separable extension $L/k$ the group $CH_0(X_L)$.  By means of duality in Voevodsky's derived category of motives, one can show (see \cite[Theorem 2.2]{HuberKahn} or \cite[\S 3.4]{Deglise2}) that for $L$ as above, there is a canonical identification
\[
\H_0^S(X)(L) \isomto CH_0(X_L).
\]
If one replaces $CH_0(X_L)$ by its rationalization, a similar statement is true for Suslin homology with $\Q$-coefficients.  The sections of $\H_0^S(X)$ over a smooth scheme $U$ with function field $L$ can thus be described either in terms of unramified elements, or as follows.  If $u$ is a codimension $1$ point of $U$, there are specialization maps $CH_0(X_L) \to CH_0(X_{\kappa_{\nu}})$ \cite[\S 20.3]{Fulton}, and $\H_0^S(X)(U)$ can be realized as the intersection of the kernels of these specialization maps.
\end{ex}

\begin{rem}
By the equivalence of categories between an appropriate category of strictly $\aone$-invariant sheaves with transfers and Rost's category of cycle modules \cite[Th\'eor\`eme 3.3]{Deglise1}, $\H_0^S(X)$ gives rise to Merkurjev's universal cycle module from \cite[\S 2.3]{Merkurjev}.  Lemma \ref{lem:h0universal} combined with this observation can be used to give an alternate proof of \cite[Theorem 2.11]{Merkurjev} under the hypothesis that $k$ is perfect.
\end{rem}

\begin{ex}
\label{ex:unramifiedwittgroup}
Let $k$ be a field having characteristic unequal to $2$.  For any smooth $k$-scheme $X$, let $W(X)$ be the associated Witt group. Using the purity results of \cite{OjangurenPanin}, one can study the Nisnevich sheafification of this presheaf.  Indeed, the Nisnevich sheafification of the functor $X \mapsto W(X)$ defines a sheaf ${\bf W}$, which we refer to as the unramified Witt sheaf.  One can identify the group of sections ${\bf W}(X)$ as the subgroup of $W(k(X))$ with trivial (second) residues at points of codimension $1$ of $X$.  See \cite[Appendice]{CTO}, \cite[\S 2.1]{MMilnor} and the references therein for more details.  While Witt groups do not admit transfers in the same sense as Milnor K-theory or unramified \'etale cohomology, there is a notion of transfer for Witt groups.
\end{ex}

\begin{ex}
\label{ex:unramifiedfundamentalideal}
Let ${\bf W}$ be the unramified Witt sheaf as just defined.  Let $I(k)$ denote the fundamental ideal in the Witt ring of $k$, i.e., the ideal of even dimensional forms, and let $I^n$ denote the $n$-th power of the fundamental ideal (which is known to be additively generated by Pfister forms.  We then set ${\bf I}^n(X) = I^n(k(X)) \cap {\bf W}(X)$.  The presheaf $U \mapsto {\bf I}^n(U)$ is a strictly $\aone$-invariant sheaf by, e.g., \cite[Theorem 2.3]{MMilnor}.  There is a monomorphism of strictly $\aone$-invariant sheaves ${\bf I}^{n+1} \hookrightarrow {\bf I}^{n}$ (coming from the corresponding injective maps on sections over fields), and it follows that ${\bf I}^n/{\bf I}^{n+1}$ is also strictly $\aone$-invariant.
\end{ex}

\begin{ex}
\label{ex:unramifiedbrauergroup}
If $k$ is a field having characteristic exponent $p$, let $\gm'$ denote the \'etale sheaf $\gm \tensor_{\Z} \Z[\frac{1}{p}]$.  Let ${\mathbf H}^2_{\et}(\gm')$ denote the Nisnevich sheaf associated with the presheaf $U \mapsto H^2_{\et}(U,\gm')$.  One can show that $H^2_{\et}(U,\gm')$ is strictly $\aone$-invariant.  Using Lemma \ref{lem:h0universal} one deduces that
\[
\hom_{\Ab^{\aone}_k}(\H_0^{\aone}(X),{\mathbf H}^2_{\et}(\gm')) = H_{Nis}^0(X,{\mathbf H}^2_{\et}(\gm')).
\]
Furthermore, one can show using purity that if $X$ is smooth and proper $H_{Nis}^0(X,{\mathbf H}^2_{\et}(\gm'))$ is precisely the cohomological Brauer group $H^2_{\et}(X,\gm')$.  It follows from Proposition \ref{prop:dependsonpi0} and Lemma \ref{lem:h0universal} that if $X$ is an $\aone$-connected smooth scheme over an algebraically closed field, then $H_{Nis}^0(X,{\mathbf H}^2_{\et}(\gm'))$ is trivial.  That $H^2_{\et}(X,\gm')$ is trivial if $X$ is $\aone$-connected was first observed by B. Bhatt; this result is stated (with proof) in \cite[Theorem 4.3]{Gille}.  For yet another proof of this statement, see \cite[Proposition 4.2]{AM}.
\end{ex}

Given any object in the stable $\aone$-homotopy category (i.e., a $\pone$-spectrum), Morel's connectivity results (recalled here as Theorem \ref{thm:stableaoneconnectivity}) show that the associated stable $\aone$-homotopy sheaves are strictly $\aone$-invariant.  Thus, given any cohomology theory representable in the stable $\aone$-homotopy category, one can get corresponding strictly $\aone$-invariant sheaves; this applies notably to motivic cohomology, algebraic K-theory, Hermitian K-theory, etc.  Another notable example comes from the stable $\aone$-homotopy groups of motivic spheres; the known computations are related to the Milnor-Witt K-theory sheaves mentioned in Example \ref{ex:aonehomologyofgm}.

\subsubsection*{Detecting $\aone$-connectedness with $\aone$-homology and birational sheaves}
Finally, we provide the ``unramified" characterization of $\aone$-connectedness stated in the introduction as Theorem \ref{thmintro:characterization}.  This result can be viewed as an extension of \cite[Theorem 1]{AH}, and the techniques are similar.

\begin{thm}
\label{thm:characterization}
If $k$ is a field, $R$ is a commutative unital ring (e.g., $\Z$ or $\Q$) and $X$ is a smooth proper $k$-scheme, then $X$ is $\aone$-connected if and only if the canonical map $\H_0^{\aone}(X,R) \to R$ is an isomorphism.
\end{thm}

\begin{proof}
We prove only the statement with $\Z$-coefficients; the corresponding statement with $R$-coefficents follows by repeating the argument word for word with $\Z$ replaced by $R$.  If $X$ is $\aone$-connected, then the canonical map in question is an isomorphism by Proposition \ref{prop:dependsonpi0} (note: this does not require properness).  In the other direction, suppose $X$ is not $\aone$-connected.  It suffices to provide a strictly $\aone$-invariant sheaf $M$ such that the map $M(k) \to M(X)$ is not an isomorphism.  By Lemma \ref{lem:h0universal} this is equivalent to proving that the map
\[
\hom_{\Ab^{\aone}_k}(\Z,M) \longrightarrow \hom_{\Ab^{\aone}_k}(\H_0^{\aone}(X),M)
\]
is not a bijection.

Recall that a presheaf of sets $\F$ on $\Sm_k$ is called birational if for any open dense immersion $U \hookrightarrow X$ the map $\F(X) \to \F(U)$ is a bijection.  In \cite[Theorem 6.2.1]{AM}, we showed that if $X$ is a smooth proper $k$-scheme, then there are a birational and $\aone$-invariant sheaf $\pi_0^{b\aone}(X)$ and a morphism $X \to \pi_0^{b\aone}(X)$ (functorial in morphisms of smooth proper schemes) characterized by the property that if $L$ is any finitely generated separable extension of $k$, then $\pi_0^{b\aone}(X)(L) = X(L)/R$.  Here, the set $X(L)/R$ is the set of $R$-equivalence classes of points in $X(L)$.  Now, either $X(k)$ is empty or not.

{\em Case 1}. Suppose $X$ is $\aone$-disconnected, but $X(k)$ is empty.  In \cite[Lemma 2.4]{AH}, we proved that the free sheaf of abelian groups on $\pi_0^{b\aone}(X)$, denoted $\Z(\pi_0^{b\aone}(X))$, is birational and strictly $\aone$-invariant.  Homomorphisms $\Z \to \Z(\pi_0^{b\aone}(X))$ correspond precisely to elements of $\pi_0^{b\aone}(X)(k)$.  In \cite[Corollary 2.9]{AH}, we showed that $X(k)$ is non-empty if and only if the map $\H_0^{\aone}(X) \to \Z$ is an epimorphism.

{\em Case 2}.  Assume that $X$ is $\aone$-disconnected, but $X(k)$ is non-empty.  Any rational point in $X(k)$ induces a splitting $\Z \to \H_0^{\aone}(X)$, and a corresponding splitting $\Z \to \Z(\pi_0^{b\aone}(X))$.  Since the category of strictly $\aone$-invariant sheaves of groups is abelian (see Corollary \ref{cor:strictlyaoneinvariantabeliancategory}), we have direct sum decompositions $\H_0^{\aone}(X) \cong \Z \oplus \tilde{\H}_0^{\aone}(X)$ and $\Z(\pi_0^{b\aone}(X)) \cong \Z \oplus \widetilde{\Z(\pi_0^{b\aone}(X))}$, and these splittings are compatible in the sense that the morphism $\H_0^{\aone}(X) \to \Z(\pi_0^{b\aone}(X))$ induces a morphism $\Z \oplus \tilde{H}_0^{\aone}(X) \to \Z \oplus \widetilde{\Z(\pi_0^{b\aone}(X))}$ that is the identity morphism on the first summand.

By \cite[Corollary 2.4.4]{AM}, $X$ is $\aone$-connected if and only if for every finitely generated separable extension $L/k$ the set $\pi_0^{b\aone}(X)(L)$ is reduced to a point.  Thus, by assumption, there exists a separable extension $K/k$ such that $\pi_0^{b\aone}(X)(K)$ consists of (strictly) more than $1$ element.  Write $X_K$ for the base extension of $X$ to $\Spec K$.  Pullback gives an identification $\H_0^{\aone}(X)(K) = \H_0^{\aone}(X_K)(K)$ by \cite[\S 5.1]{MStable}, see in particular Example 5.1.3.  Thus, without loss of generality, we can assume $k = K$ and that $\pi_0^{b\aone}(X)(k)$ consists of strictly more than $1$ element.

Each element of $\pi_0^{b\aone}(X)(k)$ determines a homomorphism $\Z \to \H_0^{\aone}(X)$ that is non-trivial, since the composite morphism $\Z \to \H_0^{\aone}(X) \to \Z(\pi_0^{b\aone}(X))$ is non-trivial.  Taking the sum of these homomorphisms gives rise to a non-trivial homomorphism $\Z(\pi_0^{b\aone}(X))(k) \to \H_0^{\aone}(X)(k)$.  It follows immediately that $\tilde{\H}_0^{\aone}(X)(k)$ is non-trivial.

\end{proof}

Combining Lemma \ref{lem:unramifiedinvariants} with Theorem \ref{thm:characterization} we deduce the following result.

\begin{cor}
If $k$ is a field and $X$ is a smooth proper $k$-scheme, then $X$ is $\aone$-connected if and only if for every strictly $\aone$-invariant sheaf $M$, the canonical map $M(k) \to M^{ur}(X)$ is a bijection.
\end{cor}

\subsubsection*{Detecting $\aone$-connectedness with Suslin homology}
As it turns out, the key point in the proof of Theorem \ref{thm:characterization} is the use of strictly $\aone$-invariant sheaves that do not necessarily possess transfers.  In general, e.g., if $X$ is a smooth curve of genus $g \geq 1$, the sheaves $\Z(\pi_0^{b\aone}(X))$ do not possess transfers; for a discussion of this point, see \cite[\S 2]{LevineST}.  If $X$ is a smooth proper $k$-variety, neither the zeroth Suslin homology sheaf of $X$ with integral nor the variant with rational coefficients can detect $\aone$-connectedness.  With rational coefficients, this fits into a general statement about Suslin homology of separably rationally connected varieties (see, e.g., \cite[Chapter 4 Definition 3.2]{Kollar}) and appeal to Example \ref{ex:artinmumford}.  If $k$ is a perfect field, \cite[Chapter 4 Theorem 3.9.4]{Kollar} shows that separably rationally connected varieties $X/k$ have the property that for every separably closed extension $L/k$, every two $L$-points can be connected by a $\pone$.

\begin{prop}
\label{prop:rationallyconnected}
If $k$ is a perfect field, and $X$ is a smooth proper $k$-scheme such that for every separably closed field $L/k$ we have $X(L)/R = \ast$, then the canonical morphism $\H_0^{S}(X,\Q) \to \Q$ is an isomorphism.
\end{prop}

\begin{proof}
Since the sheaf $\H_0^S(X,\Q)$ is strictly $\aone$-invariant, it suffices by Corollary \ref{cor:isomorphismsofstrictlyaoneinvariantsheaves} to prove that the map in question is an isomorphism on sections over every finitely generated separable extension $L/k$.  By means of the identification $\H_0^S(X,\Q)(L) = CH_0(X_L)_{\Q}$ from Example \ref{ex:unramifiedchow} it therefore suffices to prove $CH_0(X_L)_{\Q}$ is isomorphic to $\Q$ under the stated hypotheses.

If $\bar{L}$ is an algebraic closure of $L$, then we have a restriction map for Chow groups
\[
CH_0(X_L) \to CH_0(X_{\bar{L}}).
\]
The kernel of this map is a torsion subgroup.  Indeed, if $Z$ is a cycle in $CH_0(X_L)$ that goes to zero in $CH_0(X_{\bar{L}})$, then $Z$ necessarily goes to zero in a finite extension $L'/L$.  On the other hand pullback followed by pushforward is multiplication by $[L':L]$, and thus $[L':L] Z = 0$.

Under the assumption on $X$, we know that $CH_0(X_{\bar{L}}) = \Z$.  Upon tensoring with $\Q$, restriction becomes an isomorphism: it is surjective since $X_L$ has a $0$-cycle of finite degree coming from a point over some finite extension $L'/L$ and injective since the kernel of restriction is torsion and therefore becomes trivial after tensoring with $\Q$.
\end{proof}

\begin{ex}
\label{ex:artinmumford}
The classic examples of \cite{ArtinMumford} provide unirational (hence separably rationally connected) smooth proper varieties $X$ over $\cplx$ that are non-rational, but for which $Br(X)$ is non-trivial.  In particular, these varieties have $\H_0^{S}(X,\Q) = \Q$, but are not $\aone$-connected, e.g., by Example \ref{ex:unramifiedbrauergroup} and Corollary \ref{cor:aonelocaldetection}.  Other examples along these lines are provided in \cite{CTO, Peyre1} and \cite{ABK}.
\end{ex}

Finally, we observe that the zeroth integral Suslin homology sheaf of a smooth projective variety cannot detect $\aone$-connectedness.  Since we know that the zeroth $\aone$-homology detects rational points, and the zeroth Suslin homology is related to $0$-cycles, a natural place to look for a counterexample is among the smooth projective $k$-varieties that possess a $0$-cycle of degree $1$ but that have no $k$-rational point; we thank Sasha Merkurjev for pointing out the following example due to Parimala.

\begin{ex}
\label{ex:parimala}
If $X$ is a smooth projective variety such that the morphism $\H_0^S(X) \to \Z$ is an isomorphism, it need not be the case that $X$ is $\aone$-connected.  In \cite[Theorem 3]{Parimala}, Parimala gives a field $k$ (having characteristic $0$) and a projective homogeneous space $X$ under a connected reductive linear algebraic group over $k$ such that i) $X$ has a point over a degree $2$ and degree $p$ ($p$ odd) extension of $k$, but ii) has no $k$-rational point.  Point (ii) guarantees that $X$ is not $\aone$-connected.

Again combining \ref{ex:unramifiedchow} and Corollary \ref{cor:isomorphismsofstrictlyaoneinvariantsheaves}, to prove that $\H_0^S(X) \to \Z$ is an isomorphism, it suffices to prove that $CH_0(X_L) \to \Z$ is an isomorphism for every finitely generated extension, separable extensions $L/k$.  If $K$ is a field, and $X$ is a projective homogeneous space under a connected reductive linear algebraic group such that $X(K)$ is non-empty, then choice of $x \in X(K)$ determines an isomorphism $G/P \isomt X$, where $P$ is a $K$-parabolic subgroup of $G$.  Then, \cite[V Theorem 21.20(ii)]{Borel} states that $X_K$ is $K$-rational, and $K$-birational invariance of the Chow group of $0$-cycles (see Remark \ref{rem:birationalinvarianceofchow}) allows one to deduce that $CH_0(X_K) = \Z$.  Combining this discussion with point (i), it follows that for any extension $L/k$, the degree map $CH_0(X_L) \to \Z$ is also an isomorphism.  In fact, since having a $0$-cycle of degree $1$ is equivalent to having points over extensions of coprime degrees, this argument shows that if $X$ is any projective homogeneous space under a connected reductive group that has a $0$-cycle of degree $1$, then $\H_0^S(X) \to \Z$ is an isomorphism.
\end{ex}

\begin{footnotesize}
\bibliographystyle{alpha}
\bibliography{birationalinvariants}
\end{footnotesize}
\end{document}